\documentclass[tbtags,reqno]{amsart}
\usepackage{graphicx,color}
\usepackage{amsmath,amsfonts,amstext,amsbsy,amsopn,amsthm,amscd,amsxtra,dsfont,amssymb,pb-diagram}
\usepackage{amsmath,amsthm,amsopn,amsfonts,amssymb,epsfig,enumerate}
\usepackage{amsthm,amsopn,amsfonts,amssymb,epsfig,mathrsfs,cancel,wasysym}
\usepackage[all]{xy}
\usepackage{tabmac}
\usepackage[active]{srcltx}
\usepackage{tikz}
\usepackage{color}

\numberwithin{equation}{section}

\makeatletter
\renewcommand{\subsubsection}{\@startsection
{subsubsection} {3} {0mm} {\baselineskip} {-0.5\baselineskip} {\normalfont\normalsize\bfseries}} \makeatother

\newtheorem{theorem}{Theorem}
\newtheorem{lemma}[theorem]{Lemma}
\newtheorem{proposition}[theorem]{Proposition}

\newtheorem{corollary}[theorem]{Corollary}

\def\cal L{{\mathcal L}}

\def\cal L{{\mathcal L}}

\def \part {\vdash}

\def\beq{\begin{equation}}
\def\eeq{\end{equation}}

\begin{document}

\title[MultiMacdonald polynomials]
{Multi-Macdonald polynomials}

\author{Camilo Gonz\'alez}
\address{Instituto de Matem\'atica y F\'{\i}sica, Universidad de Talca, Casilla 747, Talca,
Chile} \email{cgonzalez@inst-mat.utalca.cl}

\author{Luc Lapointe}
\address{Instituto de Matem\'atica y F\'{\i}sica, Universidad de Talca, Casilla 747, Talca,
Chile} \email{lapointe@inst-mat.utalca.cl}


\begin{abstract}
  We introduce Macdonald polynomials indexed by $n$-tuples of partitions and characterized by certain orthogonality and triangularity relations.  We prove that they can be explicitly given as products of ordinary Macdonald polynomials depending on special alphabets.  With this factorization in hand, we establish their most basic properties, such as explicit formulas for their norm-squared, evaluation and reproducing kernel.  Moreover, we show that the $q,t$-Kostka coefficients associated to the multi-Macdonald polynomials are positive and correspond to $q,t$-analogs of the dimensions of the irreducible representations of $C_n \sim S_d$, the wreath product of the cyclic group $C_n$ with the symmetric group.

\end{abstract}

\keywords{Macdonald polynomials}

\maketitle

\section{Introduction}

Let $n \geq 1$ be a fixed integer.  Our goal is to define Macdonald polynomials
in the space of multi-symmetric functions in the 
sets of variables (or alphabets) $x^{(i)}=x_1^{(i)}, x_2^{(i)}, \ldots$ of infinite cardinality for $i$ from $1$ to $n$.  To be more precise, let a multi-partition be an
$n$-tuple of partitions $\pmb{\lambda} = (\lambda^{(1)}, \lambda^{(2)}, \ldots, \lambda^{(n)})$.  The space of multi-symmetric function is the $\mathbb Q(q,t)$-vector space
whose basis is given by the multi-monomial symmetric functions
\begin{equation}
m_{\pmb \lambda}(x^{(1)},\ldots ,x^{(n)}):=m_{\lambda^{(1)}}(x^{(1)}) m_{\lambda^{(2)}}(x^{(2)}) \cdots m_{\lambda^{(n)}}(x^{(n)}) \, ,
\end{equation}
where $m_\lambda(x)$ is the usual monomial symmetric function (see Section~\ref{SecDef} for the basic concepts in symmetric function theory).
Observe that the fact that the cardinality of the alphabets is infinite ensures that there are no relations among the multi-monomial symmetric functions.
The dominance ordering on multi-partitions is the following: $\pmb \mu \leq \pmb \lambda$ iff
\begin{align*}
 \left(\sum_{i < k} \big|\mu^{(i)}\big|\right) + \mu^{(k)}_1+ \cdots + \mu^{(k)}_j \leq \left(\sum_{i < k} \big|\lambda^{(i)}\big|\right) + \lambda^{(k)}_1+ \cdots + \lambda^{(k)}_j  \quad {\rm ~for~all~} j {\rm ~and~} k. 
 \end{align*}
Moreover for each multi-partition ${\pmb \lambda}=(\lambda^{(1)},\lambda^{(2)},\ldots,\lambda^{(n)})$ we associate the composition $\overline{\pmb\lambda}:=(|\lambda^{(1)}|,|\lambda^{(2)}|,\ldots,|\lambda^{(n)}|)$, while for each composition $\alpha\in \mathbb{N}^n$, we let $|\alpha|= \alpha_1+\cdots + \alpha_n$ and $n(\alpha)=\sum_i (i-1)\alpha_i$.

Before enunciating our main theorem, we need to define the multi-power sum symmetric functions
\begin{align*}
p_{\pmb \lambda}(x^{(1)},\ldots,x^{(n)}):=p_{\lambda^{(1)}}(x^{(1)}) p_{\lambda^{(2)}}(x^{(1)},x^{(2)})\cdots 
p_{\lambda^{(n)}}(x^{(1)},\ldots,x^{(n)}) \, ,
\end{align*}
where again $p_\lambda(x)$ are the usual power-sum symmetric functions.
Observe how intertwined the sets of variables are this time.

The multi-Macdonald polynomials have a triangularity/orthogonality characterization reminiscent of that of the Macdonald polynomials.
\begin{theorem}\label{main-theorem}
There exists a unique basis $\{P_{\pmb \lambda}^{(q,t)}(x^{(1)},\ldots,x^{(n)})\}_{\pmb \lambda}$ such
that
\begin{align}\label{tri-scalar}
1)& \quad  P_{\pmb \lambda}^{(q,t)}(x^{(1)},\ldots,x^{(n)})=
\,m_{\pmb \lambda}(x^{(1)},\ldots, x^{(n)}) + {\rm ~lower~terms}
\\
2) & \quad \langle \! \langle
 P_{\pmb \lambda}^{(q,t)},P_{\pmb \mu}^{(q,t)}\rangle\!\rangle_{q,t,n}=
0 \quad {\rm if} \quad  {\pmb \lambda}\neq {\pmb \mu},
\end{align}
where the scalar product is given by
\begin{align}\label{Innerproduct1}
\langle\!\langle
p_{\pmb \lambda},p_{\pmb \mu}\rangle\!\rangle_{q,t,n}=\delta_{\pmb \lambda ,
\pmb \mu}  z_{\pmb \lambda}(q,t)  := \delta_{\pmb \lambda ,
\pmb \mu}  q^{(n-1)|\overline{\pmb\lambda}|-n(\overline{\pmb\lambda})}
 z_{\lambda^{(1)}} \cdots z_{\lambda^{(n)}}
\prod_{i=1}^{\ell(\lambda^{(n)})}\frac{1-q^{\lambda_i^{(n)}}}{1-t^{\lambda_i^{(n)}}},
\end{align}
with $z_\lambda= \prod_{i \geq 1}  i^{m_i(\lambda)} m_i(\lambda)!$ if $m_i(\lambda)$ is the number of
entries equal to $i$ in the partition $\lambda$.
\end{theorem}
In order to prove the theorem, we actually construct
a basis that satisfies the two properties stated in the theorem (the uniqueness follows immediately from the uniqueness of the Gram-Schmidt process when implemented using any linear order compatible with the dominance ordering on multi-partitions).  This basis is obtained as a product of usual Macdonald polynomials albeit at very special alphabets.
\begin{proposition} \label{propfac} 
  Using the plethystic notation (see Section~\ref{SecDef}),
  the multi-Macdonald polynomials can be given explicitly as
\begin{align} \label{macdoexp}
P_{\pmb \lambda}^{(q,t)}(x^{(1)},\ldots,x^{(n)})=P_{\lambda^{(1)}}^{(q,q^{n-1}t)}\left[A_{n-1}\right] P_{\lambda^{(2)}}^{(q^{n-1}t,q^{n-2}t)}
\left[A_{n-2} \right]  \cdots    P_{\lambda^{(n)}}^{(qt,t)}\left[A_0\right]
\end{align}
where the alphabets $A_i$ are defined
recursively, starting from $A^{(0)}=X^{(n)}$, as
\begin{equation}\label{multi-recursion}
A_i=X^{(n-i)}+\frac{q(1-q^{i-1}t)}{1-q^i t}A_{i-1}
\end{equation}
for $1\leq i\leq n-1$. 
\end{proposition}  
We point out that the case $n=2$ was studied in \cite{BLM}, in which case the multi-Macdonald polynomials were called double Macdonald polynomials. 

Owing to the factorization \eqref{macdoexp}, we can establish many properties
of the multi-Macdonald polynomials.  We show that the invariance of the Macdonald polynomials when $q,t \mapsto q^{-1},t^{-1}$ has a natural extension to the multi-Macdonald case (see Proposition~\ref{symmetrytheorem}). We obtain a reproducing kernel for the scalar product \eqref{Innerproduct1}, as well as 
explicit formulas for their norm-squared and evaluation.  Furthermore, we show that the $q,t$-Kostka coefficients associated to the multi-Macdonald polynomials are positive and correspond to $q,t$-analogs of the dimensions of the irreducible representations of $C_n \sim S_d$, the wreath product of the cyclic group $C_n$ with the symmetric group.  This suggests that multi-Macdonald polynomials can be considered as wreath product Macdonald polynomials (another construction of wreath product Macdonald polynomials is presented in \cite{Hai2}).

\section{Basic definitions} \label{SecDef}

We first recall some definitions
related to partitions \cite{Mac}.
A partition $\lambda=(\lambda_1,\lambda_2,\dots)$ of degree $|\lambda|$
is a vector of non-negative integers such that
$\lambda_i \geq \lambda_{i+1}$ for $i=1,2,\dots$ and such that
$\sum_i \lambda_i=|\lambda|$.  The length $\ell(\lambda)$
of $\lambda$ is the number of non-zero entries of $\lambda$.
Each partition $\lambda$ has an associated Ferrers' diagram
with $\lambda_i$ lattice squares in the $i^{th}$ row,
from the top to bottom. Any lattice square in the Ferrers diagram
is called a cell (or simply a square), where the cell $(i,j)$ is in the $i$th row and $j$th
column of the diagram.  
The conjugate $\lambda'$ of  a partition $\lambda$ is represented  by
the diagram
obtained by reflecting  $\lambda$ about the main diagonal.
Given a cell $s=(i,j)$ in $\lambda$, we let 
\begin{equation} \label{eqarms}
  a_{\lambda}(s)=\lambda_i-j\, , \qquad l_{\lambda}(s)=\lambda_j'-i \,
\qquad a'_\lambda(s)=j-1 \qquad {\rm and} \qquad l'_\lambda(s)=i-1
  .
\end{equation}
The quantities $a_{\lambda}(s)$ and  $a'_{\lambda}(s)$ are respectively called
the arm-length and the arm-colength while $l_{\lambda}(s)$ and  $l'_{\lambda}(s)$ are respectively called
the leg-length and the leg-colength.

The dominance ordering on partitions is defined such that
$$
\mu \leq \lambda  \qquad {\rm iff} \qquad \mu_1 + \cdots +\mu_k \leq \lambda_1 + \cdots +\lambda_k \quad \forall k  
$$

Before defining the Macdonald polynomials, we need to introduce two bases of the ring of symmetric functions. For $x=x_1,x_2,\dots,x_N$, the monomial symmetric functions are such that
$$
m_\lambda(x)= \sum_\alpha x_1^{\alpha_1} \cdots x_N^{\alpha_N}
$$
where the sum is over all distinct permutations $\alpha$ of $(\lambda_1,\dots,\lambda_N)$ (if necessary, a string of 0's is added at the end of $\lambda$).  In the following we shall always consider that $N$ is infinite.  Using the $i^{th}$ power-sum 
$$
p_i(x)= x_1^i+x_2^i + \cdots
$$
the power-sum basis is then simply defined as
$$
p_\lambda(x) = \prod_{i=1}^{\ell(\lambda)} p_{\lambda_i}(x)
$$
For our purposes, the ring of symmetric functions will be simply taken as the ring $\Lambda_{\mathbb K}=\mathbb K[p_1(x),p_2(x),\dots]$ over any field $\mathbb K$.

The Macdonald polynomials 
$\{P_{ \lambda}^{(q,t)}(x)\}_{ \lambda}$ depend on two parameters $q$ and $t$ and
form the unique basis of the ring $\Lambda_{\mathbb Q(q,t)}$ (in the remainder of the article, $\mathbb K$ will always be equal to $\mathbb Q(q,t)$)  such that
\begin{align} \label{mactri}
1)& \quad  P_{\lambda}^{(q,t)}(x)=
\,m_{\lambda}(x) + {\rm ~lower~terms}
\\  
2) & \quad \langle \! \langle
 P_{\lambda}^{(q,t)},P_{\mu}^{(q,t)}\rangle\!\rangle_{q,t}=
0 \quad {\rm if} \quad  { \lambda}\neq { \mu},
\end{align}
where the scalar product is given on the power-sums by
\begin{align} \label{macscalprod}
\langle\!\langle
p_{\lambda},p_{\mu}\rangle\!\rangle_{q,t}=\delta_{\lambda ,
\mu}  z_{\lambda}(q,t)  := \delta_{\lambda ,
  \mu} z_{\lambda}
\prod_{i=1}^{\ell(\lambda)}\frac{1-q^{\lambda_i}}{1-t^{\lambda_i}},
\end{align}
with $z_\lambda$ defined in Theorem~\ref{main-theorem}.  The existence of the Macdonald poplynomials is non-trivial and follows from the construction of a difference operator $D$ (the Macdonald operator) whose eigenvectors are the Macdonald polynomials.

It will prove very convenient for our purposes to use the language of $\lambda$-rings (or plethysms).
The power-sum $p_i$ acts on the ring of rational
formal power series in $q,t,x_1,x_2, \ldots$ with coefficient in the
field $\mathbb{Q}$ as
\begin{align*}
p_i\left[ \frac{\sum_\alpha c_\alpha u_\alpha} {\sum_\beta
d_\beta v_\beta}\right]=\frac{\sum_\alpha c_\alpha
u_\alpha^i}{\sum_\beta d_\beta v_\beta^i},
\end{align*}
where $c_\alpha, d_\beta\,\,\in\,\, \mathbb{Q}$ and where
$u_\alpha, v_\beta$ are monomials in $q,t,x_1,x_2, \ldots$.  Since
the power-sums form a basis of the ring of symmetric functions, this
extends uniquely to an action of the ring of symmetric functions on
the ring of rational formal power series in $q,t,x_1, x_2, \ldots$.
In
this notation, a symmetric function $f(x)$ is denoted $f[X]$,
where $X=x_1+x_2+\cdots$.  Similarly, 
letting $X^{(i)}=x_1^{(i)}+ x_2^{(i)}+ \cdots $ for $i=1,\dots,n$,
the multi-symmetric functions $m_{\pmb \lambda}(x^{(1)},\dots,x^{(n)})$ and $p_{\pmb \lambda}(x^{(1)},\dots,x^{(n)})$
will for instance be respectively denoted
$m_{\pmb \lambda}[X^{(1)},X^{(2)},\dots ,X^{(n)}]$ and $p_{\pmb \lambda}[X^{(1)},X^{(1)}+X^{(2)},\dots ,X^{(1)}+\cdots +X^{(n)}]$.

\section{Proof of Proposition~\ref{propfac}}

As previously mentioned, the case $n=2$ was studied in \cite{BLM}.
The proof relied in this case on the following lemma that will again prove crucial in this article.
\begin{lemma}[\cite{BLM}]  \label{lemmadouble} In the case $n=2$, the scalar product $\langle\!\langle \cdot , \cdot \rangle\!\rangle_{q,t,2}$ is identical to the scalar product
$\langle  \cdot , \cdot \rangle_{q,t,2}$  defined as
\begin{align} \label{eqlemma}
  \bigl \langle & \,   p_{\pmb \lambda}\left[A_{1},A_{0}\right]
    , 
p_{\pmb \mu}\left[A_{1} ,A_0 \right]
\,  \bigr \rangle_{q,t,2}  := \delta_{\pmb \lambda,\pmb \mu}  q^{|\lambda^{(1)}|} z_{\lambda^{(1)}}(q,qt)  z_{\lambda^{(2)}}(qt,t)
\end{align}
or, more explicitly, as
\begin{align} \label{eqlemma2}
  \bigl \langle  \,   p_{\lambda^{(1)}}\left[X^{(1)}+\frac{q(1-t)}{(1-qt)} X^{(2)}\right]  p_{\lambda^{(2)}}\left[X^{(2)}\right]
   & ,  p_{\mu^{(1)}}\left[X^{(1)}+\frac{q(1-t)}{(1-qt)} X^{(2)}\right]  p_{\mu^{(2)}}\left[X^{(2)}\right]
  \,  \bigr \rangle_{q,t,2} \nonumber \\
&  \qquad \qquad := \delta_{\pmb \lambda,\pmb \mu}  q^{|\lambda^{(1)}|} z_{\lambda^{(1)}}(q,qt)  z_{\lambda^{(2)}}(qt,t)
\end{align}
\end{lemma}
The lemma has the following analog in the multi-Macdonald case.
\begin{lemma} \label{lemmaorthon}
the scalar product 
$\langle\!\langle \cdot , \cdot \rangle\!\rangle_{q,t,n}$ is equal to the scalar product defined by
\begin{align} \label{scalar2}
  \bigl \langle & \,   p_{\pmb \lambda}\left[A_{n-1},A_{n-2},\dots,A_0\right]
    , 
p_{\pmb \mu}\left[A_{n-1}, A_{n-2} ,\dots ,A_0 \right]
 \,  \bigr \rangle_{q,t,n} \nonumber \\
 & \qquad \qquad := \delta_{\pmb \lambda,\pmb \mu}  
  q^{(n-1)|\overline{\pmb\lambda}|-n(\overline{\pmb\lambda})} z_{\lambda^{(1)}}(q,q^{n-1}t)  z_{\lambda^{(2)}}(q^{n-1}t,q^{n-2}t) \cdots  z_{\lambda^{(n)}}(qt,t)
\end{align}
\end{lemma}  
\noindent {\it Proof:} 
The lemma amounts to showing that
$\langle\!\langle \cdot , \cdot \rangle\!\rangle_{q,t,n}$ is equal to the scalar product
\begin{align} \label{scalar22}
  \bigl \langle & \,   p_{\pmb \lambda}\left[A_{n-1},A_{n-2},\dots,A_0\right]
    , 
p_{\pmb \mu}\left[A_{n-1}, A_{n-2} ,\dots ,A_0 \right]
 \,  \bigr \rangle_{q,t,n} \nonumber \\
 & \qquad  :=  
 q^{(n-1)|\overline{\pmb\lambda}|-n(\overline{\pmb\lambda})}
  \bigl \langle p_{\lambda^{(1)}}[A_{n-1}]\, , \,   p_{\mu^{(1)}}[A_{n-1}]  \bigr \rangle_{q,q^{n-1}t}
\prod_{i=2}^n \,  \bigl \langle p_{\lambda^{(i)}}[A_{n-i}]\, , \,  p_{\mu^{(i)}}[A_{n-i}]  \bigr \rangle_{q^{n-i+1}t,q^{n-i}t}
 \end{align}
where by abuse of notation we always consider that $ \bigl \langle p_{\lambda}[Z]\, , \,  p_{\mu}[Z]  \bigr \rangle_{q,t} = \delta_{\lambda \mu} z_\lambda (q,t)$ for any alphabet $Z$ (which in our case corresponds to $Z=A_0,\dots,A_{n-1}$).

We will proceed by induction starting from the case $n=2$ which is covered by Lemma~\ref{lemmadouble}.  For the general case $n> 2$, we decompose the scalar product \eqref{scalar22} as
\begin{align}
   \bigl \langle  \,    p_{\pmb \lambda} & \left[A_{n-1}, A_{n-2},\dots,A_0\right]
    ,  
 p_{\pmb \mu}\left[A_{n-1}, A_{n-2} ,\dots ,A_0 \right]
 \,  \bigr \rangle_{q,t,n} = \bigl\langle \,  p_{\lambda^{(n)}}
\left[A_{0} \right] ,  p_{\mu^{(n)}}
\left[A_{0} \right] \, 
\bigr \rangle_{qt,t}    \nonumber \\
  & \qquad  \times 
q^{|\lambda^{(1)}|+\cdots+|\lambda^{(n-1)}|}
 \bigl \langle  \,   p_{\pmb {\hat \lambda}}\left[A_{n-1},A_{n-2},\dots,A_1\right]
    , 
p_{\pmb {\hat \mu}}\left[A_{n-1}, A_{n-2} ,\dots ,A_1 \right]
 \,  \bigr \rangle_{q,qt,n-1}
\end{align}
where $\pmb {\hat \lambda}=(\lambda^{(1)},\dots,\lambda^{(n-1)})$.
Since the $q$-power is constant on multi-partitions $\pmb {\hat \lambda}$ of the same total degree, we can use induction to write
\begin{align}
\bigl \langle  \,    p_{\pmb \lambda} & \left[A_{n-1}, A_{n-2},\dots,A_0\right]
    ,  
 p_{\pmb \mu}\left[A_{n-1}, A_{n-2} ,\dots ,A_0 \right]
 \,  \bigr \rangle_{q,t,n} = \bigl\langle \,  p_{\lambda^{(n)}}
\left[A_{0} \right] ,  p_{\mu^{(n)}}
\left[A_{0} \right] \, 
\bigr\rangle_{qt,t} \nonumber \\
& \quad \times q^{|\lambda^{(1)}|+\cdots+|\lambda^{(n-1)}|}
\langle\!\langle  p_{\pmb {\hat \lambda}}\left[ Y^{(1)}, \dots, Y^{(n-2)},Z \right]
,  p_{\pmb {\hat \mu}}\left[ Y^{(1)}, \dots, Y^{(n-2)},Z \right]
  \rangle\!\rangle_{q,qt,n-1}
\end{align}
where $Y^{(i)}=X^{(1)}+\cdots + X^{(i)}$ and $Z= Y^{(n-1)} + q(1-t)X^{(n)}/(1-qt)=Y^{(n-2)}+A_1$.
By the expression for the scalar product \eqref{Innerproduct1}, the part that depends on the alphabets
$A_0=X^{(n)}$ and $Z$ and on the partitions $\lambda^{(n-1)}$ and  $\lambda^{(n)}$ 
yields
\begin{align}
q^{|\lambda^{(n-1)}|} & \bigl\langle \,  p_{\lambda^{(n-1)}}
\left[Z \right] ,  p_{\mu^{(n-1)}}
\left[Z \right] \, 
\bigr\rangle_{q,qt} \bigl\langle \,  p_{\lambda^{(n)}}
\left[X^{(n)} \right] ,  p_{\mu^{(n)}}
\left[X^{(n)} \right] \, 
\bigr\rangle_{qt,t} \nonumber \\
& \qquad = \langle \! \langle \, p_{\lambda^{(n-1)}}[Y^{(n-1)}] p_{\lambda^{(n)}}[Y^{(n)}],  p_{\mu^{(n-1)}}[Y^{(n-1)}] p_{\mu^{(n)}}[Y^{(n)}] \, \rangle \! \rangle_{q,t,2}
\end{align}
from Lemma~\ref{lemmadouble} (observe that we used the fact that the $q$-power appearing in  $\langle\!\langle  \cdot
,  \cdot
  \rangle\!\rangle_{q,qt,n-1}$ does not depend on $\lambda^{(n-1)}$ and  $\lambda^{(n)}$). Hence
\begin{align}
\bigl \langle  \,    p_{\pmb \lambda} & \left[A_{n-1}, A_{n-2},\dots,A_0\right]
    ,  
 p_{\pmb \mu}\left[A_{n-1}, A_{n-2} ,\dots ,A_0 \right]
 \,  \bigr \rangle_{q,t,n}
 \nonumber
 \\
  & \quad
  =  q^{(n-1)|\overline{\pmb\lambda}|-n(\overline{\pmb\lambda})} z_{\lambda^{(1)}}  z_{\lambda^{(2)}} \cdots  z_{\lambda^{(n-1)}}  z_{\lambda^{(n)}}(q,t)
 \\
  & \quad
  = \langle\!\langle  p_{\pmb {\lambda}}\left[ Y^{(1)}, \dots, Y^{(n-1)},Y^{(n)} \right]
,  p_{\pmb {\mu}}\left[ Y^{(1)}, \dots, Y^{(n-1)}, Y^{(n)} \right]
  \rangle\!\rangle_{q,t,n}
\end{align}
which proves the lemma.
\begin{flushright}
$\square$
\end{flushright}

\noindent {\it Proof of Proposition~\ref{propfac}:}
We have to show that the expression given in \eqref{macdoexp}
satisfies the triangularity and the orthogonality in Theorem~\ref{main-theorem}.
The orthogonality is an immediate consequence of Lemma~\ref{lemmaorthon}
since the scalar product is a product of usual Macdonald scalar products \eqref{macscalprod} at the right alphabets (the extra $q$-power does not affect the
orthogonality since it only depends on the degree of the components of the multi-partition).

We now consider the triangularity.
Suppose by induction that the triangularity in the monomial basis holds in the case when there are $n-1$ sets of alphabets (the base case $n=1$ is the usual Macdonald case), i.e.
\begin{align}
 P_{\pmb {\tilde\lambda}}^{(q,t)}
=\sum_{\pmb {\tilde \mu} \leq \pmb {\tilde \lambda}} *   \,\,m_{\pmb {\tilde \mu}}(x^{(2)},\ldots,x^{(n)})
\end{align}
where we use the tilde in $\pmb {\tilde \lambda}$ to emphasize that the multi-partition
has $n-1$ components instead of $n$, and where, for simplicity, we use $*$ to denote the expansion coefficients (instead of for instance the more cumbersome $c_{\pmb {\tilde \lambda}, \pmb {\tilde \mu}}(q,t)$).

For the case with $n$ alphabets, we need to analyse the extra factor $P_{\lambda^{(1)}}^{(q^nt,q^{n-1}t)}\left[A_{n-1}\right]$.  It is triangular, from \eqref{mactri}, in the monomial basis   
\begin{align*}
P_{\lambda^{(1)}}^{(q^nt,q^{n-1}t)}\left[A_{n-1}\right]=\sum_{\nu^{(1)}\leq \lambda^{(1)}} *\,\,m_{\nu^{(1)}}\left[ A_{n-1}\right]
\end{align*}
We now have to expand $m_{\nu^{(1)}}\left[ A_{n-1}\right]$ in the $m_{\pmb \mu}(x^{(1)},\ldots,x^{(n)})$ basis.  Using the notation  $(\mu,\nu) \subseteq \lambda$ to denote $\mu \subseteq \lambda$ and $\nu \subseteq \lambda$ (the same notation will be used in the general case with more than two components), we have 
\begin{align*}
m_{\nu^{(1)}}\left[ X^{(1)}+\frac{q-q^{n-1}t}{1-q^{n-1}t}A_{n-2}\right]&=\sum_{ |(\rho^{(1)},\gamma)|=|\nu^{(1)}| \, : \, (\rho^{(1)},\gamma)\subseteq \nu^{(1)}} m_{\rho^{(1)}}[X_1] m_{\gamma}\left[\frac{q-q^{n-1}t}{1-q^{n-1}t}A_{n-2} \right]\\
&=\sum_{ |(\rho^{(1)},\sigma)|=|\nu^{(1)}| \, :\, \rho^{(1)}\subseteq \nu^{(1)}}*\,\, m_{\rho^{(1)}}[X^{(1)}] m_{\sigma}\left[A_{n-2} \right]
\end{align*}
where $|(\mu,\nu)|=|\mu|+|\nu|$.  Repeating this argument again and again, we obtain
\begin{align*}
m_{\nu^{(1)}}\left[ A_{n-1}\right]=
\sum_{|\pmb \rho|=|\nu^{(1)}| \, :  \, \pmb {\hat \rho}\subseteq \nu^{(1)}} \, m_{\pmb {\rho}}(x^{(1)},\ldots,x^{(n)})
\end{align*}
where, as before, $\pmb {\hat \rho}= (\rho^{(1)},\dots,\rho^{(n-1)})$.
Hence, we have
\begin{align*}
  P_{\pmb \lambda}^{(q,t)}(x^{(1)},\ldots,x^{(n)})
=\sum_{\nu^{(1)}\leq \lambda^{(1)}}\,\,\, \sum_{\pmb {\tilde \mu}\leq \pmb {\tilde \lambda}}\,\,\,\,\,\, \sum_{|\pmb \rho|=|\nu^{(1)}| \, :  \, \pmb {\hat \rho}\subseteq \nu^{(1)}}
  * \, m_{\pmb \rho}(x^{(1)},\ldots,x^{(n)}) \, m_{\pmb {\tilde \mu}}(x^{(2)},\ldots,x^{(n)})
\end{align*}
It is known that 
\begin{align*}
m_\mu[X] \, m_{\nu}[X]= m_{\mu+\nu}[X]+\sum_{\sigma<\mu+\nu} \,*\,\, m_\sigma[X],
\end{align*}
where $\mu+\nu=(\mu_1+\nu_1,\mu_2+\nu_2,\dots)$.  Therefore
\begin{align*}
   P_{\pmb \lambda}^{(q,t)}(x^{(1)},\ldots,x^{(n)})
   =\sum_{\nu^{(1)}\leq \lambda^{(1)}}\,\,\, \sum_{\pmb {\tilde \mu}\leq \pmb {\tilde \lambda}}\,\,
   \,\, \sum_{|\pmb \rho|=|\nu^{(1)}| \, :  \, \pmb {\hat \rho}\subseteq \nu^{(1)}}
\, \, \sum_{\pmb {\tilde \sigma} \preceq \pmb {\tilde \rho} + \pmb {\tilde \mu}, \sigma^{(1)}=\rho^{(1)}}
   * \, m_{\pmb {\sigma}}(x^{(1)},\ldots,x^{(n)})
\end{align*}
where $\pmb {\tilde \sigma} \preceq \pmb {\tilde \rho} + \pmb {\tilde \mu}$ stands for $\sigma^{(2)} \leq \rho^{(2)}+\mu^{(2)}, \ldots, \sigma^{(n)} \leq \rho^{(n)}+\mu^{(n)}$.

We will now see that the conditions on the summation indices in the previous equation imply that $\pmb \sigma \leq \pmb \lambda$. From
$\sigma^{(1)}=\rho^{(1)}$, $\rho^{(1)} \subseteq \nu^{(1)}$ and $\nu^{(1)} \leq \lambda^{(1)}$, we first deduce that 
\begin{align}
\sum_{i=1}^k \sigma_i^{(1)} = \sum_{i=1}^k \rho_i^{(1)} \leq  \sum_{i=1}^k \nu_i^{(1)} \leq  \sum_{i=1}^k \lambda_i^{(1)} 
\end{align}  
We then let $2 \leq \ell \leq n$ and use  $\sigma^{(1)}=\rho^{(1)}$ as well as $\pmb {\tilde \sigma} \preceq \pmb {\tilde \rho} + \pmb {\tilde \mu}$ to obtain
\begin{align}
  |\sigma^{(1)}|  + \cdots & + |\sigma^{(\ell-1)}| + \sum_{i=1}^k \sigma^{(\ell)}_i \nonumber \\
  & \leq
  |\rho^{(1)}| +\left(|\rho^{(2)}|+|\mu^{(2)}| \right) +\cdots + \left( |\rho^{(\ell-1)}|
  + |\mu^{(\ell-1)}| \right)+ \sum_{i=1}^k \left( \rho^{(\ell)}_i+ \mu^{(\ell)}_i \right) 
\end{align}
Then, using $|\pmb \rho|=|\nu^{(1)}|$ followed by $\pmb {\tilde \mu}\leq \pmb {\tilde \lambda}$ and $|\nu^{(1)}|=|\lambda^{(1)}|$, we get that
\begin{align}
  |\sigma^{(1)}| + \cdots + |\sigma^{(\ell-1)}|+ \sum_{i=1}^k \sigma^{(\ell)}_i & 
\leq
  |\nu^{(1)}| +|\mu^{(2)}| +\cdots + |\mu^{(\ell-1)}|+ \sum_{i=1}^k \mu^{(\ell)}_i \nonumber \\
  & \leq
  |\lambda^{(1)}| +|\lambda^{(2)}| +\cdots + |\lambda^{(\ell-1)}|+ \sum_{i=1}^k \lambda^{(\ell)}_i  
\end{align}  
which proves the triangularity.
\begin{flushright}
$\square$
\end{flushright}

\section{Properties of multi-Macdonald polynomials}

We now establish the properties of the Macdonald polynomials that extend to the  multi-Macdonald case.

\subsection{Sending $q \mapsto q^{-1}$ and $t \mapsto t^{-1}$.}
The Macdonald polynomials safisfy the following property \cite{Mac}
\begin{align}\label{symmetrymacdonald}
P_\lambda^{(q^{-1},t^{-1})} (x)= P_\lambda^{(q,t)} (x).
\end{align}
The generalization of that property to the multi-Macdonald case is the following.
\begin{proposition}\label{symmetrytheorem} Using the plethystic notation, the multi-Macdonald polynomials are such such that
  \begin{align*}
    P_{\pmb \lambda}^{(q^{-1},t^{-1})}[X^{(1)},X^{(2)},\ldots,X^{(n)}]= q^{n(\overline{\pmb\lambda})-(n-1)|\overline{\pmb\lambda}|} 
    P_{\pmb \lambda}^{(q,t)} [q^{n-1}X^{(1)},q^{n-2}X^{(2)},\ldots,X^{(n)}]
    \end{align*}
\end{proposition}

If $f$ is a multi-symmetric function, we define the homomorphism $\varphi$ as
$$\varphi(f[X^{(1)},\ldots,X^{(i)},\ldots,X^{(n)})]=f[q^{n-1}X^{(1)},\ldots,q^{n-i}X^{(i)},\ldots,X^{(n)}].$$
In order to prove the proposition, we first prove the following lemma:
\begin{lemma}\label{Lemmasymmetry}
  The alphabets $A_i:=A_i(q,t)$ defined recursively in \eqref{multi-recursion}
  are such that
\begin{align}
  \varphi\bigl(A_i(q,t) \bigr)=q^i A_i(q^{-1},t^{-1}) \, .
\end{align}
\end{lemma}
\noindent{\it Proof:}\,\, We proceed induction.  The result holds trivially in the base case $A_0$ given that $A_0=X^{(n)}$.
Supposing that 
\begin{align*}
   \varphi\bigl(A_i(q,t) \bigr)=q^i A_i(q^{-1},t^{-1})
\end{align*}
we then get
\begin{align*}
  \varphi\bigl(A_{i+1}(q,t)\bigr)& =q^{i+1}X^{(n-(i+1))}+\frac{q(1-q^{i}t)}{1-q^{i+1}t} \cdot q^i A_i(q^{-1},t^{-1}) \\
  & = q^{i+1} \Bigl(X^{(n-(i+1))}+\frac{q^{-1}(1-q^{-i}t^{-1})}{1-q^{-i-1}t^{-1}} \cdot  A_i(q^{-1},t^{-1}) \Bigr) = q^{i+1} A_{i+1}(q^{-1},t^{-1})
\end{align*}
which proves the lemma.
\begin{flushright}
$\square$
\end{flushright}

\noindent{\it Proof of Proposition~\ref{symmetrytheorem}:}\,\,
Using  \eqref{symmetrymacdonald}, the explicit form \eqref{macdoexp} of the multi-Macdonald polynomial and the previous lemma, we have
\begin{align*}
&  P_{\pmb \lambda}^{(q^{-1},t^{-1})}[X^{(1)},X^{(2)},\ldots,X^{(n)}] \\
&\quad  = P_{\lambda^{(1)}}^{(q,q^{n-1}t)}\left[A_{n-1}(q^{-1},t^{-1})\right] P_{\lambda^{(2)}}^{(q^{n-1}t,q^{n-2}t)}
  \left[A_{n-2}(q^{-1},t^{-1}) \right]  \cdots    P_{\lambda^{(n)}}^{(qt,t)}\left[A_0(q^{-1},t^{-1})\right] \\
  &\quad  = P_{\lambda^{(1)}}^{(q,q^{n-1}t)}\left[q^{-n+1}\varphi(A_{n-1})\right] P_{\lambda^{(2)}}^{(q^{n-1}t,q^{n-2}t)}
  \left[q^{-n+2}\varphi(A_{n-2}) \right]  \cdots    P_{\lambda^{(n)}}^{(qt,t)}\left[\varphi(A_0)\right] 
\end{align*}
Using the homogeneity of the Macdonald polynomials (which implies that $P_\lambda^{(q,t)}[q^aX]=q^{a|\lambda|} P_\lambda^{(q,t)}[X] $) and the fact that $\varphi$ is a homomorphism, we then get
\begin{align*}
    P_{\pmb \lambda}^{(q^{-1},t^{-1})}[X^{(1)},X^{(2)},\ldots,X^{(n)}]= q^{n(\overline{\pmb\lambda})-(n-1)|\overline{\pmb\lambda}|} 
    \varphi \bigl( P_{\pmb \lambda}^{(q,t)} [X^{(1)},X^{(2)},\ldots,X^{(n)}] \bigr)
    \end{align*}
This proves the proposition.
\begin{flushright}
$\square$
\end{flushright}

\subsection{Norm}
 As we will see, the explicit form \eqref{macdoexp} of a multi-Macdonald polynomial leads to an explicit expression for its norm-squared $||P_{\pmb \lambda}^{(q,t)}||^2_{q,t,n}:= \langle\!\langle
P_{\pmb \lambda}^{(q,t)} , P_{\pmb \lambda}^{(q,t)} \rangle\!\rangle_{q,t,n}$.  First, recall that Macdonald polynomials are such that \cite{Mac}
\begin{equation}
 \langle\!\langle
 P_{\lambda}^{(q,t)} , P_{\lambda}^{(q,t)} \rangle\!\rangle_{q,t} = b_\lambda(q,t)
\end{equation}
 where
\begin{align*}
b_\lambda(q,t)=\prod_{s\in \lambda}\frac{1-q^{a(s)}t^{l(s)+1}}{1-q^{a(s)+1}t^{l(s)}}.
\end{align*}
with $a(s)=a_\lambda(s)$ and $l(s)=l_\lambda(s)$ such as defined in Section~\ref{SecDef}.

Using \eqref{macdoexp}, the explicit formula for the norm-squared of the Macdonald polynomials implies from Lemma~\ref{lemmaorthon} a similar expression for the norm-squared of the multi-Macdonald polynomials.
\begin{corollary} \label{coronorm}
The Multi-Macdonald polynomial $P_{\pmb \lambda}^{(q,t)}(x^{(1)},\ldots,x^{(n)})$ is such that:
\begin{align*}
||P_{\pmb \lambda}^{(q,t)}||^2_{q,t,n}=q^{(n-1)|\overline{\pmb \lambda}|-n(\overline{\pmb \lambda})}b_{\lambda^{(1)}}(q,q^{n-1}t)^{-1} b_{\lambda^{(2)}}(q^{n-1}t,q^{n-2}t)^{-1} \cdots b_{\lambda^{(n)}}(qt,t)^{-1}
\end{align*}
\end{corollary}

\subsection{Kernel} The Macdonald polynomial scalar product has the following reproducing kernel
\begin{equation}
K(x,y;q,t) = \prod_{i,j }\frac{(t x_i y_j;q)_\infty}{(x_i y_j;q)_\infty}
\end{equation}  
where $x$ (resp. $y$) stands for $x_1,x_2,x_3,...$ (resp. $y_1,y_2,y_3,...$)  and
where
$$
(a;q)_{\infty}= \prod_{i=0}^\infty (1-q^{i}a)
$$
Being a reproducing kernel, $K(x,y;q,t)$ is such that
\begin{equation} \label{eqkern}
K(x,y;q,t) = \sum_{\lambda} \frac{1}{z_{\lambda}(q,t)} p_\lambda(x) p_\lambda(y)
\end{equation}

We now extend this result to the multi-Macdonald case.
For each $k\in \{1,\ldots,n\}$, let $\pmb x^{(k)}$ (resp. $\pmb y^{(k)}$) be the union of the alphabets $x^{(1)}, x^{(2)},\ldots, x^{(k)}$ (resp. $y^{(1)},y^{(2)},\ldots,y^{(k)}$); to simplify the notation, when $k=n$ we write $\pmb x$ (resp. $\pmb y$) instead of $\pmb x^{(n)}$(resp. $\pmb y^{(n)}$).  Even though, the alphabets  $x^{(1)}, x^{(2)},\ldots, x^{(k)}$ are infinite, the alphabet $\pmb x^{(k)}$ is countably infinite and we will suppose that its elements are ordered as $\pmb x^{(k)}_1,\pmb x^{(k)}_2,\dots$ (the order is irrelevant). Now, let
\begin{align*}
K(\pmb x, \pmb y;q,t)=\prod_{i,j} \frac{(t \pmb x_i \pmb y_j ;q)_\infty}{(\pmb x_i, \pmb y_j ;q)_\infty}\prod_{k=1}^{n-1}\prod_{i,j}\frac{1}{1-q^{-(n-k)}\pmb x^{(k)}_i \pmb y^{(k)}_j}
\end{align*}
\begin{proposition} We have that
  $$
  K(\pmb x, \pmb y;q,t) = \sum_{\pmb \lambda} \frac{1}{z_{\pmb \lambda}(q,t)}
   p_{\pmb \lambda}(x^{(1)},\ldots,x^{(n)}) p_{\pmb \lambda}(y^{(1)},\ldots,y^{(n)})
   $$
where $ z_{\pmb \lambda }(q,t)$ was defined in \eqref{Innerproduct1}. 
      As such, 
$K(\pmb x, \pmb y;q,t)$ is a reproducing kernel for the scalar product \eqref{Innerproduct1}, which implies in particular
\begin{equation}
K(\pmb x, \pmb y;q,t)= \sum_{\pmb \lambda} \frac{1}{||P_{\pmb \lambda}^{(q,t)}||^2_{q,t,n}} P_{\pmb \lambda}^{(q,t)}(x^{(1)},\ldots,x^{(n)}) P_{\pmb \lambda}^{(q,t)}(y^{(1)},\ldots,y^{(n)})
  \end{equation}
where $||P_{\pmb \lambda}^{(q,t)}||^2_{q,t,n}$ is given explicitly in Corollary~\ref{coronorm}.
\end{proposition}
\noindent {\it Proof:}
Recall that the Cauchy identity is such that  \cite{Mac}
$$
\prod_{i,j} \frac{1}{1-x_i y_j} = \sum_{\lambda} \frac{1}{z_\lambda} p_\lambda(x) p_\lambda (y)
$$
By homogeneity, the Cauchy identity is easily seen to be such that
$$
\prod_{i,j} \frac{1}{1-u x_i y_j} = \sum_{\lambda} \frac{1}{z_\lambda} u^{|\lambda|} p_\lambda(x) p_\lambda (y)
$$
The proposition then follows from \eqref{eqkern} since
\begin{align*}
K(\pmb x, \pmb y;q,t)=&\left(\sum_{\lambda^{(n)}} z_{\lambda^{(n)}}(q,t)^{-1}p_{\lambda^{(n)}}(\pmb x) p_{ \lambda^{(n)}}(\pmb y)\right) 
 \left( \sum_{\lambda^{(n-1)}}z_{\lambda^{(n-1)}}^{-1}q^{-|\lambda^{(n-1)}|} p_{\lambda^{(n-1)}}(\pmb x^{(n-1)})  p_{\lambda^{(n-1)}}(\pmb y^{(n-1)})  \right)\\
&\qquad \qquad \times \cdots \times \left( \sum_{\lambda^{(1)}}z_{\lambda^{(1)}}^{-1}q^{-(n-1)|\lambda^{(1)}|} p_{\lambda^{(1)}}(\pmb x^{(1)})  p_{\lambda^{(1)}}(\pmb y^{(1)})  \right)\\
=&\sum_{\pmb \lambda} z_{\pmb \lambda }(q,t)^{-1} p_{\pmb \lambda}(x^{(1)},\ldots,x^{(n)}) p_{\pmb \lambda}(y^{(1)},\ldots,y^{(n)}) 
\end{align*}
\begin{flushright}
$\square$
\end{flushright}

\subsection{Specializations}
We now describe the various specializations of the multi-Macdonald polynomials presented in the figure below.
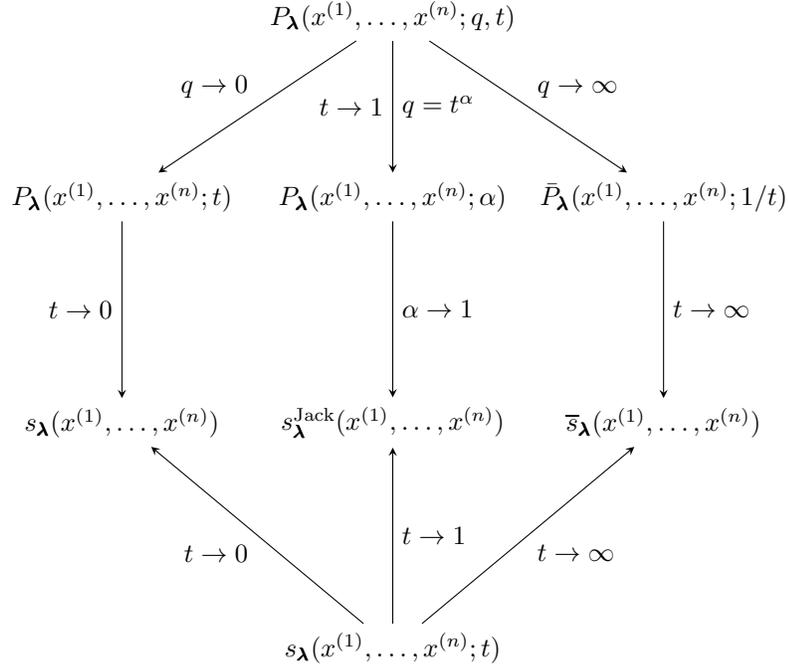
\begin{center}
  \begin{tikzpicture}
      [scale=.6,auto,>=stealth]
  \node (n1) at (2,16)	{$P_{\pmb \lambda}(x^{(1)},\ldots,x^{(n)};q,t)$};
  \node (n2) at (-4,12)	{$P_{\pmb \lambda}(x^{(1)},\ldots,x^{(n)};t)$};
  \node (n3) at (2,12)	{$P_{\pmb \lambda}(x^{(1)},\ldots,x^{(n)};\alpha)$};
  \node (n4) at (8,12)	{$\bar P_{\pmb \lambda}(x^{(1)},\ldots,x^{(n)};1/t)$};
  \node (n5) at (-4,7)	{$s_{\pmb \lambda}(x^{(1)},\ldots,x^{(n)})$};
  \node (n6) at (2,7)	{$s^{\text{Jack}}_{\pmb\lambda}(x^{(1)},\ldots,x^{(n)})$};
  \node (n7) at (8,7)	{$\overline{s}_{\pmb\lambda} (x^{(1)},\ldots,x^{(n)})$};
  \node (n8) at (2,2)	 {$s_{\pmb \lambda}(x^{(1)},\ldots,x^{(n)};t)$};
    
    \draw[->]  (n1) to node [swap]{$q\rightarrow 0$} (n2);
    \draw[->]  (n1) to node  {$q=t^\alpha$} node [swap] {$t\rightarrow 1$} (n3);
    \draw[->]  (n1) to node  {$q\rightarrow \infty$} (n4);
    \draw[->]  (n2) to node [swap] {$t\rightarrow 0$}  (n5);
    \draw[->]  (n3) to node  {$\alpha\rightarrow 1$}  (n6);
    \draw[->]  (n4) to node  {$t\rightarrow \infty$}  (n7);
    \draw[->]  (n8) to node {$t\rightarrow 0$} (n5);
    \draw[->]  (n8) to node [swap]{$t\rightarrow 1$} (n6);
    \draw[->]  (n8) to node [swap]{$t\rightarrow \infty$} (n7);
\end{tikzpicture}
\end{center}

\noindent {\it Multi-Jack polynomials:} The Multi-Jack polynomias are simply
products of usual Jack polynomials $P_\lambda^{(\alpha)}$
by taking the limit $q=t^\alpha$, $t\rightarrow 1$ in \eqref{macdoexp}.  To be more explicit,
we let the alphabets $B_i$ be defined recursively as
\begin{align*}
B_i=X^{(n-i)}+\frac{\alpha(i-1)+1}{\alpha i+1}B_{i-1}
\end{align*}
(starting from $B_0=X^{(n)}$).  We then get
\begin{align*}
P_{\pmb \lambda}^{(\alpha)}(x^{(1)},\ldots,x^{(n)})=P^{(\alpha_{n-1})}_{\lambda^{(1)}}[B_{n-1}] P_{\lambda^{(2)}}^{(\alpha_{n-2})}[B_{n-2}]\cdots P_{\lambda^{(n)}}^{(\alpha_{0})}[B_0]
\end{align*}
where, for $k=0,\dots,n-1$, we have
$$
\alpha_k =
\left \{
\begin{array}{ll}
  \displaystyle{\frac{\alpha}{\alpha(n-1)+1}} & {\rm if~} k=n-1 \\
  & \\
   \displaystyle{\frac{\alpha (k+1) +1}{\alpha k+1} }& {\rm otherwise}
\end{array}
\right .
$$

\noindent{\it Multi-Hall Littlewood polynomials: }
It is known \cite{Mac} that $P_\lambda^{(0,0)}(x)=s_\lambda(x)$ is a Schur function and that
$P_\lambda^{(0,t)}(x)=P_\lambda(x;t)$ is a Hall-Littlewood polynomial.  Now, 
when letting $q=0$ in $\eqref{multi-recursion}$, the recursion trivializes and we get $A_i=X^{(n-i)}$.  Hence, we get from \eqref{macdoexp} that
\begin{align*}
P_{\pmb \lambda}(x^{(1)},\ldots,x^{(n)};t)=s_{\lambda^{(1)}}(x^{(1)}) s_{\lambda^{(2)}}(x^{(2)})\cdots s_{\lambda^{(n-1)}}(x^{(n-1)}) P_{\lambda^{(n)}} (x^{(n)};t),
\end{align*}

When $q\rightarrow \infty$, it is straightforward to obtain that the alphabets
$A_i$ in \eqref{multi-recursion} become
$$
C_{i}= \lim_{q\to \infty} A_{n-i} = X^{(i)}+ X^{(i+1)}+\cdots+ X^{(n-1)}+(1-1/t)X^{(n)}
$$
for $i=1,\dots,n-1$, while $A_0$ remains equal to $X^{(n)}$.  Hence, we have
\begin{align*}
\bar P_{\pmb \lambda}(x^{(1)},\ldots,x^{(n)};1/t)=&s_{\lambda^{(1)}}[C_{1}] \, s_{\lambda^{(2)}}[C_2] \cdots \,  s_{\lambda^{(n-1)}}[C_{n-1}]P_{\lambda^{(n)}}(x^{(n)}; 1/t)
\end{align*}
since $P_\lambda^{(\infty,t)}(x)=P_\lambda(x;1/t)$.

\noindent{\it Multi-Schur functions: } When setting $q=t$,
the dependency over $t$ does not disappear (contrary to the Macdonald case)
and we obtain a family of multi-Schur functions $s_{\pmb \lambda}(x^{(1)},\ldots,x^{(n)};t)$ depending on $t$.

Setting $q=t$ in $\eqref{multi-recursion}$, the recursion becomes
$$
A_i=X^{(n-i)}+\frac{t(1-t^i)}{1-t^{i+1}} A_{i-1}
$$
which implies that
\begin{equation} \label{Aischur}
A_i = \sum_{k=0}^i \frac{t^{i-k}(1-t^{k+1})}{1-t^{i+1}} X^{(n-k)}
\end{equation}
in this case.  Letting $D_i=A_i$ as above, we have from \eqref{macdoexp} that
\begin{equation} \label{schurt}
s_{\pmb \lambda}(x^{(1)},\ldots,x^{(n)};t)= P_{\lambda^{(1)}}^{(t,t^n)}[D_{n-1}] P_{\lambda^{(2)}}^{(t^n,t^{n-1})}[D_{n-2}] \cdots P_{\lambda^{(n)}}^{(t^2,t)}[D_{0}]  
\end{equation}
We stress that the specializations of the Macdonald polynomials appearing in the previous product have not, to the best of our knowledge, been studied and do not appear to have any special properties.

We now describe in more details the specializations $t=1$, $t=0$ and $t=\infty$
of the multi-Schur functions. Taking the limit $t\to 1$ in \eqref{Aischur}, and setting
$$
E_i = \sum_{k=0}^i \frac{{k+1}}{{i+1}} X^{(n-k)}
$$
we obtain from \eqref{schurt} that
$$
s^{\rm Jack}_{\pmb \lambda}(x^{(1)},\ldots,x^{(n)})= P_{\lambda^{(1)}}^{(1/n)}[E_{n-1}] P_{\lambda^{(2)}}^{(n/n-1)}[E_{n-2}] \cdots P_{\lambda^{(n-1)}}^{(3/2)}[E_{1}] P_{\lambda^{(n)}}^{(2)}[E_{0}]
$$
Note that this is a product of Jack polynomials at special values of $\alpha$.

Taking the limit $t\to 0$ in \eqref{Aischur} and \eqref{schurt}, we get
\begin{equation} \label{schur}
s_{\pmb \lambda}(x^{(1)},\ldots,x^{(n)})= s_{\lambda^{(1)}}(x^{(1)}) s_{\lambda^{(2)}}(x^{(2)})\cdots s_{\lambda^{(n)}}(x^{(n)})
\end{equation}
These multi-Schur functions will be used later in Subsection~\ref{subseckostka} when studying the multi $q,t$-Kostka coefficients.

Finally, taking the limit $t\to \infty$ in \eqref{Aischur} and \eqref{schurt}, we obtain
\begin{equation*}
  \bar s_{\pmb \lambda}(x^{(1)},\ldots,x^{(n)})= s_{\lambda^{(1)}}[X^{(1)}+\cdots+X^{(n)}] s_{\lambda^{(2)}}[X^{(2)}+\cdots+X^{(n)}]\cdots s_{\lambda^{(n)}}[X^{(n)}]
\end{equation*}

\subsection{Evaluation}

Recall \cite{Mac} that the evaluation of a symmetric polynomials is a homomorphism  defined on the power-sum symmetric function $p_r$ as 
\begin{align*}
\varepsilon_{u,t}(p_r)=\frac{1-u^r}{1-t^r},
\end{align*}
where $u$ is an indeterminate. In the plethystic notation, it simply corresponds to  $p_r$ acting on the alphabet $(1-u)/(1-t)$:
\begin{align*}
p_r\left[\frac{1-u}{1-t}\right]=\frac{1-u^r}{1-t^r}. 
\end{align*}
It is known that the evaluation of a Macdonald polynomial is given by 
\begin{align}\label{EvalMac}
P_\lambda^{(q,t)}\left[\frac{1-u}{1-t}\right]=\prod_{s\in \lambda}\frac{t^{l'(s)}-q^{a'(s)}u}{1-q^{a(s)}t^{l(s)+1}}=:w_\lambda(u;q,t).
\end{align}
For multi-symmetric polynomials, we define the evaluation by
\begin{align*}
E_{{\pmb u},q,t}(m_{\pmb \lambda}\bigl(x^{(1)},\ldots,x^{(n)})\bigr)=m_{\pmb \lambda}[\mathcal X^{(1)},\ldots,\mathcal X^{(n)}]
\end{align*}
where, for the indeterminates $u_1,\dots, u_n$, we have
$$
\mathcal X^{(i)}= q^{n-i} u_1 \cdots u_{i-1}  \left( \frac{1-u_i}{1-q^{n-i}t} \right)
$$
and where it is understood that
\begin{equation} \label{multiplet}
p_r \left [ \mathcal X^{(i)} \right]=  q^{(n-i)r} u_1^r \cdots u_{i-1}^r  \left( \frac{1-u_i^r}{1-q^{(n-i)r}t^r} \right)
\end{equation}

\begin{proposition}\label{Multi-Evaluation2}
The evaluation of the multi-Macdonald polynomial $P^{(q,t)}_{\pmb \lambda}$ is given by
\begin{align*}
E_{{\pmb u};q,t} (P^{(q,t)}_{\pmb \lambda})=\left[ \prod_{i=1}^n \left({q^{n-i}u_1\cdots u_{i-1}}\right)^{|\lambda^{(i)}|} \right]  w_{\lambda^{(1)}}(v_1;q,q^{n-1}t) w_{\lambda^{(2)}}(v_2;q^{n-1}t,q^{n-2}t) \cdots w_{\lambda^{(n)}}(v_n;qt,t)
\end{align*}
where $v_i=u_i \cdots u_n$ for $i\in \{1,\ldots,n\}$.
\end{proposition}
\noindent{\it Proof:}
From \eqref{EvalMac} and \eqref{multiplet}, we have
$$
P_\lambda^{(Q,q^{n-i}t)}\left[{q^{n-i}u_1\cdots u_{i-1}}\, \frac{1-u_{i} \cdots u_n}{1-q^{n-i} t} \right]= \left({q^{n-i}u_1\cdots u_{i-1}}\right)^{|\lambda|} w_\lambda(u_i \cdots u_n;Q,q^{n-i}t)
$$
From \eqref{EvalMac} and \eqref{macdoexp}, the proposition is  then an immediate consequence of the lemma that follows.
\begin{flushright}
$\square$
\end{flushright}
\begin{lemma}  The linear map $\rho : X^{(i)} \mapsto
\mathcal X^{(i)}$, $i=1,\dots,n$,  has the following action on the $A^{(i)}$'s defined in \eqref{multi-recursion}:
\begin{align}\label{Evalua2} 
\rho (A_i)={q^{i}u_1\cdots u_{n-i-1}}\, \frac{1-u_{n-i} \cdots u_n}{1-q^{i} t}
\end{align}
for $i=0,\dots,n-1$.
\end{lemma}
\noindent {\it Proof:} We use induction on $i$. The base case $A_0=X^{(n)}$ holds since
\begin{align*}
\rho (A_0)= \mathcal X^{(n)}=u_1\cdots u_{n-1} \frac{1-u_n}{1-t}
\end{align*}
which coincides with the $i=0$ case of \eqref{Evalua2}.
Now, supposing that \eqref{Evalua2} holds, we have
\begin{align*}
 \rho (A_{i+1}) & = \mathcal X^{(n-i-1)}+ \frac{q(1-q^i t)}{(1-q^{i+1}t)} \rho(A_{i}) \\
  &=  q^{i+1} u_1 \cdots u_{n-i-2}  \left( \frac{1-u_{n-i-1}}{1-q^{i+1}t} \right)
  + \frac{q(1-q^i t)}{(1-q^{i+1}t)} {q^{i}u_1\cdots u_{n-i-1}}\, \frac{1-u_{n-i} \cdots u_n}{1-q^{i} t} \\
  &= q^{i+1}  \frac{u_1 \cdots u_{n-i-2}}{1-q^{i+1}t} \bigl( (1- u_{n-i-1})+u_{n-i-1}(1-  u_{n-i}\cdots u_{n} )\bigr)
\end{align*}
which proves the lemma by induction.
\begin{flushright}
$\square$
\end{flushright}

\subsection{Multi-Kostka coefficients} \label{subseckostka}
The integral form of the Macdonald polynomial is defined as
\begin{equation}
J_\lambda^{(q,t)}(x)= c_\lambda(q,t) P_\lambda^{(q,t)}(x)
\end{equation}  
where
\begin{equation}
c_\lambda(q,t) = \prod_{s \in \lambda} (1-q^{a(s)}t^{l(s)+1})
\end{equation}
with again $a(s)=a_\lambda(s)$ and $l(s)=l_\lambda(s)$ such as defined in Section~\ref{SecDef}.  We will also need to introduce the modified Macdonald polynomials
\begin{equation} \label{defH}
H_\lambda^{(q,t)}(x) = J_\lambda^{(q,t)}\left[ \frac{X}{1-t}\right] 
\end{equation}  
whose expansion in terms of Schur functions gives the $q,t$-Kostka polynomials
\begin{equation}
H_\lambda^{(q,t)}(x) = \sum_{\mu} K_{\mu \lambda}(q,t) \, s_\mu(x)
\end{equation}  
It is known  \cite{GH,Hai} that $K_{\mu \lambda}(q,t) \in \mathbb N [q,t]$.

Defining the integral form of the Multi-Macdonald polynomials as
\begin{align*}
J_{\pmb \lambda}^{(q,t)}(x^{(1)},\ldots,x^{(n)})=c_{\lambda^{(1)}}(q,q^{n-1}t) c_{\lambda^{(2)}}(q^{n-1}t,q^{n-2}t)\cdots c_{\lambda^{(n)}}(qt,t) P_{\pmb \lambda} (x^{(1)},\ldots,x^{(n)};q,t)
\end{align*}
we then let the modified multi-Macdonald polynomials be such that
\begin{align} \label{defmH}
H_{\pmb\lambda}^{(q,t)}(x^{(1)},\ldots,x^{(n)})=\phi\bigl(J^{(q,t)}_{\pmb \lambda}(x^{(1)},\ldots,x^{(n)})\bigr)
\end{align}
where $\phi$, which generalizes the plethystic substitution found in \eqref{defH},  has the following action on the alphabets:
\begin{align*}
X^{(i)} \mapsto X^{(i)}, \quad i =1,\dots,n-1, \qquad \qquad X^{(1)}+\cdots+ X^{(n)} \mapsto \frac{X^{(1)}+ \cdots +X^{(n)}}{1-t}
\end{align*}
\begin{lemma} We have
\begin{equation} \label{eqAK}
\phi(A_i) = \frac{q^i t \left( X^{(1)}+\cdots + X^{(n-i-1)} \right)+X^{(n-i)}+q X^{(n-i+1)} +\cdots +q^i X^{(n)} }{1-q^i t}
\end{equation}     
\end{lemma}  
\noindent {\it Proof:}  We proceed by induction starting from the case $A_0$.  We have
\begin{align*}
  \phi(A_0)=\phi(X^{(n)}) & = \phi \bigl(X^{(1)}+\cdots + X^{(n)}-(X^{(1)}+\cdots + X^{(n-1)}) \bigr) \\
  & = \frac{X^{(1)}+\cdots + X^{(n)}}{1-t} - (X^{(1)}+\cdots + X^{(n-1)}) \\
   & = \frac{t(X^{(1)}+\cdots + X^{(n-1)})+ X^{(n)}}{1-t}
\end{align*}  
which shows that the lemma holds in that case.  Then, assuming that \eqref{eqAK} holds, we obtain
\begin{align*}
  \phi(A_{i+1})& = X^{(n-i-1)}+ \frac{q(1-q^i t)}{(1-q^{i+1}t)} \phi(A_{i}) \\
  &=X^{(n-i-1)}+\frac{q^{i+1} t \left( X^{(1)}+\cdots + X^{(n-i-1)} \right)+qX^{(n-i)}+q^2 X^{(n-i+1)} +\cdots +q^{i+1} X^{(n)} }{1-q^{i+1} t} \\
  & = \frac{q^{i+1} t \left( X^{(1)}+\cdots + X^{(n-i-2)} \right)+X^{(n-i-1)}+q X^{(n-i)} +\cdots +q^{i+1} X^{(n)} }{1-q^{i+1} t}
\end{align*}  
which proves the lemma by induction.
\begin{flushright}
$\square$
\end{flushright}

The following proposition is then immediate from \eqref{defH}, \eqref{defmH} and
\eqref{macdoexp}.
\begin{proposition} \label{theoH}
  Letting
$$  
Z_i = {q^i t \left( X^{(1)}+\cdots + X^{(n-i-1)} \right)+X^{(n-i)}+q X^{(n-i+1)} +\cdots +q^i X^{(n)} }
$$
we have  
\begin{align*}
  H_{\pmb \lambda}^{(q,t)}(x^{(1)},\ldots,x^{(n)})=
  H_{\lambda^{(1)}}^{(q,q^{n-1}t)}[Z_{n-1}] \, H_{\lambda^{(2)}}^{(q^{n-1}t,q^{n-2}t)}[Z_{n-2}]
 \cdots      H_{\lambda^{(n)}}^{(qt,t)}[Z_0]
\end{align*}
\end{proposition}
We now define the multi $q,t$-Kostka coefficients $K_{\pmb \mu \pmb \lambda}(q,t)$ as
\begin{equation} \label{mkostka}
H_{\pmb \lambda}^{(q,t)} (x^{(1)},\dots,x^{(n)}) = \sum_{\pmb \mu } K_{\pmb \mu \pmb \lambda}(q,t) \, s_{\pmb \mu}(x^{(1)},\dots,x^{(n)})
\end{equation}
where the Schur function $s_{\pmb \mu}(x^{(1)},\dots,x^{(n)})
= s_{\mu^{(1)}}(x^{(1)}) \cdots  s_{\mu^{(n)}}(x^{(n)})$ was defined in \eqref{schur}.  As in the proof of Proposition~12 in \cite{BLM}, the positivity of the usual q,t-Kostka coefficients together with the Littlewood-Richardson rule implies that the multi $q,t$-Kostka coefficients are positive.
\begin{corollary}  The  multi $q,t$-Kostka coefficients are polynomials in $q$ and $t$ with nonnegative integer coefficients, that is,
$K_{\pmb \mu \pmb \lambda}(q,t) \in \mathbb N[q,t]$ for every ${\pmb \mu}$ and   
${\pmb \lambda}$.
  \end{corollary}  
Owing to $H_{\lambda}^{(1,1)}(x)=\bigl( p_1(x) \bigr)^{|\lambda|}$, we have from Proposition~\ref{theoH} that
$$
H_\lambda^{(1,1)} (x^{(1)},\dots,x^{(n)}) = \bigl(p_1[X^{(1)}+\cdots+X^{(n)}] \bigr)^{|\lambda^{(1)}|+\cdots + |\lambda^{(n)}|} = \bigl(p_1[X^{(1)}]+\cdots+p_1[X^{(n)}] \bigr)^{|\lambda^{(1)}|+\cdots + |\lambda^{(n)}|}$$
Hence \cite{Mac},
$$
H_\lambda^{(1,1)} (x^{(1)},\dots,x^{(n)})
=
\sum_{\pmb \mu } \chi^{\pmb \mu}_{\rm Id} \, s_{\pmb \mu}(x^{(1)},\dots,x^{(n)}) 
$$
where $\chi^{\pmb \mu}_{\rm Id}$ is the character of $C_n\sim S_d$ (with $d$ the total degree $|\mu^{(1)}|+\cdots + |\mu^{(n)}|$ of the multipartition $\pmb \mu$) indexed by the irreducible representation $\pmb \mu$ evaluated at the conjugacy class of the identity.  From \eqref{mkostka}, it has the following consequence.
\begin{corollary} We have that $K_{\pmb \mu \pmb \lambda}(1,1)=\chi^{\pmb \mu}_{\rm Id}$, that is, $K_{\pmb \mu \pmb \lambda}(1,1)$ is the dimension of the irreducible representation of $C_n\sim S_d$ indexed by the multi-partition $\pmb \mu$ of total degree $d$. 
  \end{corollary}

\end{document}